\documentclass[preprint,12pt]{elsarticle}
\usepackage[T1]{fontenc}
\usepackage{amsmath,amssymb,amsthm}
\usepackage{newtxtext,newtxmath}
\usepackage[a4paper,margin=28mm]{geometry}
\usepackage{microtype}
\usepackage{enumitem}
\usepackage{needspace}
\usepackage{hyperref}
\hypersetup{hidelinks,pdftitle={Convex combinations of strict triangular norms},pdfauthor={}}
\biboptions{sort&compress}
\journal{Fuzzy Sets and Systems}
\newtheorem{theorem}{Theorem}[section]
\newtheorem{proposition}[theorem]{Proposition}
\newtheorem{lemma}[theorem]{Lemma}
\newtheorem{corollary}[theorem]{Corollary}
\theoremstyle{definition}

\newtheorem{example}[theorem]{Example}

\theoremstyle{remark}

\newcommand{\AC}{\operatorname{AC}}
\newcommand{\ACloc}{\operatorname{AC}_{\mathrm{loc}}}
\newcommand{\dd}{\,\mathrm d}

\newcommand{\TP}{T_{\mathrm P}}

\newcommand{\norm}[1]{\lVert #1\rVert}
\setlist[enumerate]{label=\textup{(\roman*)},itemsep=2pt,topsep=5pt}
\allowdisplaybreaks[2]
\begin{document}
\begin{frontmatter}
\title{Convex combinations of strict triangular norms}
\author[addr1]{Yuxu Chen}
\ead{chenyuxu@scu.edu.cn}
\affiliation[addr1]{organization={School of Mathematics, Sichuan University},
city={Chengdu},
country={China}}

\begin{abstract}
We prove that a pointwise convex combination of two strict triangular norms with positive weights is associative if and only if the two norms coincide.  Together with the known ordinal-sum reduction and nilpotent case, the theorem yields the same conclusion for all continuous triangular norms and, in particular, for associative bivariate copulas. It also gives a corresponding result for representable uninorms with a common neutral element.
\end{abstract}
\end{frontmatter}

\section{Introduction}\label{sec:intro}
Triangular norms arise in fuzzy logic, probabilistic metric spaces and copula theory \cite{KMP2000}. A pointwise convex combination preserves commutativity, monotonicity and the neutral element, but need not preserve associativity. Alsina, Frank and Schweizer asked whether a convex combination of two distinct triangular norms with positive weights can again be a triangular norm \cite[Problem~5]{AFS2003}. Their Problem~6 asks whether the product copula admits such a decomposition into two associative copulas.

Discontinuous examples are known \cite{Jenei2006,DS2008}. For continuous triangular norms, Mesiar and Mesiarov\'a-Zem\'ankov\'a \cite{MM2008} exclude distinct summands with the same diagonal. Ouyang and Fang \cite{OF2008} reduce the problem to common Archimedean ordinal-sum components and exclude a convex combination of a strict and a nilpotent norm; see also \cite[p.~2 of the expanded version]{NPS2009}. Petr\'ik and Sarkoci \cite{PS2009} settle the nilpotent case. The remaining case is therefore that of two strict triangular norms.

Petr\'ik's Conjecture~1 \cite{Petrik2010} concerns this case under assumptions on boundary derivatives and absolute continuity of canonical generators; a related restricted assertion appears in \cite[Corollary~1]{NPS2009}. The unrestricted strict question is stated in \cite[Section~5]{Navara2019} and \cite[Section~5]{NP2020}. Our main result removes these additional assumptions.

\begin{theorem}\label{thm:full-strict}
Let \(T_1,T_2\) be strict triangular norms and let \(0<\lambda<1\). Then
\[
T_\lambda(x,y)=\lambda T_1(x,y)+(1-\lambda)T_2(x,y)
\]
is associative if and only if \(T_1=T_2\).
\end{theorem}

The known reductions give Corollary~\ref{cor:continuous} for all continuous triangular norms, answering \cite[Problem~1]{Navara2015}. Corollary~\ref{cor:copulas} answers Problem~6 of \cite{AFS2003} and Open Question~1 of \cite{Seliga2026}. The result for representable uninorms with a common neutral element follows from \cite[Proposition~33]{MZ2021}.

The proof first treats the case in which the generators and their inverses are locally absolutely continuous, using identities for section derivatives. For arbitrary strict norms, we construct an auxiliary triple with this regularity and the same convex identity. The construction uses measures on common levels and comparisons of section measures. The first step then applies to the auxiliary triple.

Section~\ref{sec:prelim} collects the definitions and standard analytic facts. Section~\ref{sec:strict} proves the absolutely continuous case. Sections~\ref{sec:common}--\ref{sec:regularization} construct the auxiliary triple and complete the proof. Section~\ref{sec:scope} gives the consequences and examples.

\section{Preliminaries}\label{sec:prelim}
A \emph{triangular norm}, or \emph{\(t\)-norm}, is a map \(T:[0,1]^2\to[0,1]\) that is commutative, associative and nondecreasing in each variable and satisfies \(T(x,1)=x\). These axioms imply
\begin{equation}\label{eq:bounds}
T(x,0)=0,\qquad T(x,y)\leq\min\{x,y\}.
\end{equation}
A \(t\)-norm is called \emph{strict} if it is continuous and, for each \(y>0\), the function \(x\mapsto T(x,y)\) is strictly increasing. A continuous \(t\)-norm is \emph{Archimedean} if \(T(x,x)<x\) for \(0<x<1\). It is then either strict or \emph{nilpotent}; in the latter case, for every \(x\in(0,1)\), some finite iterate of \(x\) is zero \cite{KMP2000}.

For \(y\in(0,1)\), write
\[
L_y^T(x)=T(x,y),\qquad 0\leq x\leq1.
\]
For a strict \(T\), this is an increasing homeomorphism from \([0,1]\) onto \([0,y]\).

Every strict \(t\)-norm has an additive generator, that is, a continuous strictly decreasing bijection
\[
f:(0,1]\longrightarrow[0,\infty)
\]
with \(f(1)=0\), \(f(0+)=\infty\), and
\begin{equation}\label{eq:generator}
T(x,y)=\psi(f(x)+f(y)),\qquad \psi=f^{-1},\qquad x,y>0.
\end{equation}
The values on the zero boundary are determined by \eqref{eq:bounds}. Two additive generators define the same strict \(t\)-norm if and only if they differ by a positive multiplicative constant \cite{KMP2000}. We use this freedom to normalize generators in Section~\ref{sec:strict}.

A \emph{multiplicative generator} is an increasing homeomorphism \(h:[0,1]\to[0,1]\) such that
\[
T(x,y)=h^{-1}(h(x)h(y)).
\]
The two representations are related by \(h=e^{-f}\). Multiplicative generators of a fixed strict \(t\)-norm differ by a positive power \cite{KMP2000}.

Throughout, \(m\) denotes Lebesgue measure, and ``almost everywhere'' refers to this measure. We write \(\AC([a,b])\) for the absolutely continuous functions on a compact interval and \(\ACloc(J)\) for functions absolutely continuous on every compact subinterval of an open interval~\(J\). The symbol \(\partial_i\) denotes differentiation with respect to the \(i\)th variable, and \(\norm{\cdot}_1\) denotes the norm in \(L^1(0,1)\) unless another domain is specified.

For an increasing homeomorphism \(h:[0,1]\to[0,1]\), \(\mathrm dh\) denotes its Lebesgue--Stieltjes measure, characterized by
\[
\mathrm dh((a,b])=h(b)-h(a).
\]
It has no endpoint atoms. More generally, for any increasing continuous function \(g\) on a compact interval, \(\mathrm dg\) denotes the measure of its increments, with no endpoint atoms. If \(h\) is an increasing homeomorphism and \(B\) is Borel, then
\begin{equation}\label{eq:stieltjes-image}
\mathrm dh(B)=m(h(B)),\qquad
\mathrm d(g\circ h)(B)=\mathrm dg(h(B)).
\end{equation}
Both equalities follow first on half-open intervals and then on Borel sets by uniqueness of measures. All null sets may be replaced by Borel null supersets.

For positive Borel measures \(\alpha,\beta\), the notation \(\alpha\ll\beta\) means
\[
\beta(E)=0\quad\Longrightarrow\quad\alpha(E)=0
\qquad\text{for every Borel set }E.
\]
Thus \(m\ll\rho\) means that every \(\rho\)-null set is Lebesgue-null. For a continuous increasing function \(h\), one has \(\mathrm dh\ll m\) exactly when \(h\) is absolutely continuous \cite[Proposition~3.32]{Folland1999}.

We use the Lebesgue decomposition and Tonelli theorems in the following forms.

\begin{proposition}\label{prop:lebesgue-decomposition}
Every finite positive Borel measure \(\rho\) on a compact interval \(I\) has a unique decomposition
\[
\rho=r\,m+\rho_{\mathrm s},\qquad r\in L^1(I),\quad r\geq0,
\]
where \(r\) is determined almost everywhere and \(\rho_{\mathrm s}\) is concentrated on a Borel set of Lebesgue measure zero \cite[Theorem~3.8]{Folland1999}. If \(\rho=\mathrm dh\) for a continuous increasing \(h\), then \(r=h'\) almost everywhere \cite[Proposition~3.30]{Folland1999}.
\end{proposition}

\begin{proposition}\label{prop:tonelli}
Let \(\alpha,\beta\) be finite positive Borel measures on intervals \(I,J\), and let \(F:I\times J\to[0,\infty]\) be Borel measurable. Then
\[
\int_I\left(\int_J F(x,y)\,\beta(\mathrm dy)\right)\alpha(\mathrm dx)
=\int_J\left(\int_I F(x,y)\,\alpha(\mathrm dx)\right)\beta(\mathrm dy),
\]
with the value \(+\infty\) allowed. For a signed function integrable with respect to \(\alpha\times\beta\), the same equality holds by Fubini's theorem \cite[Theorem~2.37]{Folland1999}.
\end{proposition}

The following standard facts about absolute continuity are recorded for reference. The null-set criterion is the monotone case of \cite[Theorem~3.12]{Leoni2009}; the inverse criterion is \cite[Exercise~3.21]{Leoni2009}. Composition and extension to the endpoints follow from \cite[Exercise~3.51, Corollary~3.42 and Exercise~3.7(ii)]{Leoni2009}.

\begin{proposition}\label{prop:ACfacts}
Let \(h\) be an increasing homeomorphism between compact intervals.
\begin{enumerate}
\item The function \(h\) is absolutely continuous if and only if it maps Lebesgue null sets to Lebesgue null sets.
\item If \(h\) is absolutely continuous, then \(h^{-1}\) is absolutely continuous if and only if \(h'>0\) almost everywhere.
\item The composition of two monotone absolutely continuous functions is absolutely continuous whenever defined. A continuous monotone function on a compact interval that is locally absolutely continuous in its interior is absolutely continuous on the whole interval.
\end{enumerate}
The corresponding statements hold for decreasing homeomorphisms, with the sign in~\textup{(ii)} reversed, and locally on open intervals.
\end{proposition}

In particular, for a strict generator~\(f\),
\begin{equation}\label{eq:genregular}
f\in\ACloc((0,1)),\qquad
f^{-1}\in\ACloc((0,\infty))
\end{equation}
is equivalent to
\begin{equation}\label{eq:genregular2}
f\in\ACloc((0,1)),\qquad f'(x)<0\quad\text{almost everywhere}.
\end{equation}
Neither condition imposes a finite derivative at \(0\) or~\(1\).

\section{The absolutely continuous case}\label{sec:strict}
Throughout Sections~\ref{sec:strict}--\ref{sec:regularization}, the weights satisfy \(\lambda,\mu>0\) and \(\lambda+\mu=1\). We consider a strict triple with
\begin{equation}\label{eq:convex-combination}
T_3=\lambda T_1+\mu T_2.
\end{equation}
A convex combination of strict triangular norms is itself strict whenever it is associative, since the remaining axioms, continuity and strict increase of the positive sections are preserved.

In this section all additive generators and their inverses satisfy \eqref{eq:genregular}. Section derivatives are represented by Borel difference-quotient limits, with value zero where no finite derivative exists. All derivative identities hold almost everywhere in the interiors.

\Needspace{15\baselineskip}
\begin{lemma}\label{lem:normalization}
Let \(T_1,T_2,T_3\) satisfy \eqref{eq:convex-combination}, and suppose that all three additive generators and their inverses satisfy \eqref{eq:genregular}. Then either \(T_1=T_2=T_3\), or their additive generators can be chosen so that, with
\[
v_i=-\frac1{f_i'}\quad\text{almost everywhere},
\]
one has
\begin{equation}\label{eq:vsum}
v_3=\lambda v_1+\mu v_2\quad\text{almost everywhere},
\end{equation}
and there is a sequence \(z_n\uparrow1\) satisfying
\begin{equation}\label{eq:ratios}
\frac{f_1(z_n)}{f_3(z_n)}\longrightarrow1,
\qquad
\frac{f_2(z_n)}{f_3(z_n)}\longrightarrow1.
\end{equation}
\end{lemma}
\begin{proof}
Let \(E\subset(0,1)\) be a common full-measure set on which
\(-\infty<f_i'(x)<0\) for \(i=1,2,3\). The inverse derivative formula gives, for \(x\in E\),
\begin{equation}\label{eq:expansion}
x-T_i(x,y)=v_i(x)f_i(y)+o(f_i(y))\qquad(y\uparrow1).
\end{equation}
Set \(w_1=\lambda,w_2=\mu\) and fix \(x_0\in E\). For \(i=1,2\),
\[
0\leq w_i\bigl(x_0-T_i(x_0,y)\bigr)
\leq x_0-T_3(x_0,y).
\]
Thus \eqref{eq:expansion} yields, for all \(y<1\) sufficiently close to~\(1\),
\[
0\leq\frac{f_i(y)}{f_3(y)}
\leq\frac{4v_3(x_0)}{w_i v_i(x_0)}.
\]
Choose \(z_n\uparrow1\) with
\(f_1(z_n)/f_3(z_n)\to A\geq0\) and
\(f_2(z_n)/f_3(z_n)\to B\geq0\).
Divide the deficit identity at \((x,z_n)\) by \(f_3(z_n)\). For every \(x\in E\), \eqref{eq:expansion} gives
\begin{equation}\label{eq:vunnormalized}
v_3=\lambda A v_1+\mu Bv_2\quad\text{almost everywhere}.
\end{equation}
Since \(v_3>0\) on \(E\), \((A,B)\ne(0,0)\).

If \(A=0\), then \(B>0\) and
\(f_3'=f_2'/(\mu B)\) almost everywhere. For \(0<x<t<1\), the fundamental theorem of calculus gives
\(f_3(t)-f_3(x)=(f_2(t)-f_2(x))/(\mu B)\). Letting \(t\uparrow1\) yields \(f_3=f_2/(\mu B)\), hence
\(T_3=T_2\) and \(\lambda(T_1-T_2)=0\). The case \(B=0\) follows by interchanging the indices.

If \(A,B>0\), replace \((f_1,f_2)\) by \((f_1/A,f_2/B)\). Then
\[
(v_1,v_2)\mapsto(Av_1,Bv_2),\qquad
\frac{f_1(z_n)}{A f_3(z_n)}\to1,\quad
\frac{f_2(z_n)}{B f_3(z_n)}\to1.
\]
Equations~\eqref{eq:vsum} and~\eqref{eq:ratios} follow in the new notation.
\end{proof}

\begin{lemma}\label{lem:L1}
Let \(T\) be a strict \(t\)-norm whose additive generator and inverse satisfy \eqref{eq:genregular}. For \(x\in(0,1)\) let \(p_x\) be the almost-everywhere derivative of the section \(y\mapsto T(x,y)\). Then
\begin{equation}\label{eq:positiveintegral}
p_x>0\quad\text{almost everywhere},\qquad
\int_0^1p_x(y)\dd y=x,
\end{equation}
and the map \(x\mapsto p_x\) is continuous from \((0,1)\) to \(L^1(0,1)\).
\end{lemma}
\begin{proof}
For fixed \(x\in(0,1)\), the generator formulas
\[
F_x(y)=f^{-1}(f(x)+f(y)),\qquad
F_x^{-1}(u)=f^{-1}(f(u)-f(x))
\]
and Proposition~\ref{prop:ACfacts} show that the section \(F_x=T(x,\cdot)\) and its inverse are absolutely continuous, including their endpoints. This proves \eqref{eq:positiveintegral}.

Put \(h=-(f^{-1})'\). Then \(h>0\) almost everywhere and \(\int_0^\infty h=1\). The chain rule and change of variables for monotone absolutely continuous functions \cite[Corollaries~3.50 and~3.57]{Leoni2009}, first on compact subintervals and then by exhaustion, give
\begin{equation}\label{eq:pformula}
p_x(y)=h(f(x)+f(y))(-f'(y))\quad\text{almost everywhere}
\end{equation}
and
\begin{equation}\label{eq:L1translation}
\norm{p_x-p_{x'}}_1
=\int_0^\infty |h(f(x)+s)-h(f(x')+s)|\dd s.
\end{equation}
Extend \(h\) by zero to \(\mathbb R\). Translation continuity in \(L^1\) \cite[Proposition~8.5]{Folland1999} makes the last integral tend to zero as \(x'\to x\).
\end{proof}

The normalized derivatives satisfy the following two identities.

\begin{lemma}\label{lem:compatibility}
Assume the nontrivial alternative of Lemma~\ref{lem:normalization}, and use its normalized generators. Put
\[
r=\frac{v_2}{v_1}>0\quad\text{almost everywhere},\qquad
p_x=\partial_2T_1(x,\cdot),\qquad q_x=\partial_2T_2(x,\cdot).
\]
Choose finite positive measurable values of \(r\) on its null exceptional set. Then
\begin{equation}\label{eq:compatibility}
(r(x)-r(y))\bigl(q_x(y)-r(x)p_x(y)\bigr)=0
\quad\text{for almost every }(x,y)\in(0,1)^2.
\end{equation}
Moreover, if \(r\) takes only two distinct values \(a,b>0\) almost everywhere, then
\begin{equation}\label{eq:twovaluerestriction}
(a-1)(b-1)=0.
\end{equation}
\end{lemma}
\begin{proof}
For each \(i\), let \(N_i\subset(0,1)\) be null and contain all points where \(f_i\) has no finite strictly negative derivative. Every positive section has an absolutely continuous inverse, so
\[
m\bigl((L_y^{T_i})^{-1}(N_i)\bigr)=0
\qquad(0<y<1).
\]
The chain rule and Fubini's theorem therefore give
\begin{equation}\label{eq:derivbalance}
v_i(x)\partial_1T_i(x,y)=v_i(y)\partial_2T_i(x,y)
\quad\text{almost everywhere in }(x,y).
\end{equation}
Substitute \eqref{eq:vsum} and the derivatives of \eqref{eq:convex-combination} into the identity for \(i=3\). With \(p=p_x(y),q=q_x(y)\), the result is
\[
\begin{aligned}
0={}&(\lambda+\mu r(x))
 \left(\lambda p+\mu\frac{r(y)}{r(x)}q\right)
 -(\lambda+\mu r(y))(\lambda p+\mu q)\\
={}&\frac{\lambda\mu}{r(x)}(r(x)-r(y))(r(x)p-q).
\end{aligned}
\]
Since \(\lambda\mu/r(x)>0\), this proves \eqref{eq:compatibility}.

Suppose \(r\in\{a,b\}\) almost everywhere and set
\(D=(\lambda+\mu a)(\lambda+\mu b)\). For \(r=a,b\),
\[
\frac1{\lambda+\mu r}
=\frac\lambda D+\frac{\mu ab}{Dr},\qquad
f_2'=\frac{f_1'}r,\quad
f_3'=\frac{f_1'}{\lambda+\mu r}.
\]
Thus \(f_3'=(\lambda f_1'+\mu ab f_2')/D\) almost everywhere.
Integrating from \(x\) to \(t<1\) and letting \(t\uparrow1\) gives
\begin{equation}\label{eq:twovaluegens}
f_3=\frac\lambda D f_1+\frac{\mu ab}D f_2.
\end{equation}
Evaluate at \(z_n\), divide by \(f_3(z_n)\), and use \eqref{eq:ratios}:
\[
D=\lambda+\mu ab,\qquad
0=D-\lambda-\mu ab=-\lambda\mu(a-1)(b-1).
\]
This proves \eqref{eq:twovaluerestriction}.
\end{proof}

\begin{lemma}\label{lem:ratio}
In the setting of Lemma~\ref{lem:compatibility}, the function \(r\) is constant almost everywhere.
\end{lemma}
\begin{proof}
By Fubini's theorem, there is a full-measure set \(E\subset(0,1)\)
such that, for each \(x\in E\),
\begin{equation}\label{eq:fiberrelation}
q_x(y)=r(x)p_x(y)
\quad\text{for almost every }y\text{ with }r(y)\ne r(x).
\end{equation}
Assume that \(r\) is not essentially constant.

First suppose that \(r\in\{a,b\}\) almost everywhere, with \(a<b\).
Put \(A=\{r=a\}\), \(B=\{r=b\}\), and replace \(E\) by
\(E\cap(A\cup B)\). Both \(A\) and \(B\) have positive measure.
The closures of \(E\cap A\) and \(E\cap B\), taken relative to
\((0,1)\), cover this connected interval, so they meet at some
\(x_0\in(0,1)\). Taking sequences from each set to \(x_0\),
Lemma~\ref{lem:L1} and \eqref{eq:fiberrelation} give
\[
q_{x_0}=a p_{x_0}\quad\text{a.e. on }B,
\qquad
q_{x_0}=b p_{x_0}\quad\text{a.e. on }A.
\]
Set \(P_A=\int_Ap_{x_0}\dd m>0\) and
\(P_B=\int_Bp_{x_0}\dd m>0\). By \eqref{eq:positiveintegral},
\[
P_A+P_B=x_0=bP_A+aP_B,
\qquad
1=\frac{bP_A+aP_B}{P_A+P_B}\in(a,b),
\]
contrary to \eqref{eq:twovaluerestriction}.

Otherwise, choose three pairwise disjoint Borel intervals
\(J_1,J_2,J_3\subset(0,\infty)\) such that
\(A_j=\{r\in J_j\}\) has positive measure for each \(j\).
Define
\[
R_j(x)=\frac{\int_{A_j}q_x\dd m}{\int_{A_j}p_x\dd m}
\qquad(j=1,2,3).
\]
Lemma~\ref{lem:L1} makes these functions continuous and strictly
positive. For \(x\in E\), at most one \(J_j\) contains \(r(x)\),
so \eqref{eq:fiberrelation} gives \(R_j(x)=r(x)\) for at least
two indices. Their median,
\[
R=R_1+R_2+R_3-\min_jR_j-\max_jR_j,
\]
is therefore a continuous positive function with \(R|_E=r|_E\).

Fix \(x\in(0,1)\) and choose \(x_n\in E\) tending to \(x\).
By Lemma~\ref{lem:L1}, a common subsequence satisfies
\(p_{x_n}\to p_x\) and \(q_{x_n}\to q_x\) almost everywhere.
Since \(r(x_n)=R(x_n)\to R(x)\) and \(r=R\) almost everywhere,
\eqref{eq:fiberrelation} yields
\begin{equation}\label{eq:extendedfiber}
q_x(y)=R(x)p_x(y)
\quad\text{for almost every }y\text{ with }R(y)\ne R(x).
\end{equation}
If \(R\) were nonconstant, its range would contain an interval.
Only countably many level sets have positive measure; hence there
would be a \(c\ne1\) in its range with
\(m(R^{-1}(\{c\}))=0\). Taking \(x_0\) with \(R(x_0)=c\),
\eqref{eq:extendedfiber} and \eqref{eq:positiveintegral} would give
\[
x_0=\int_0^1q_{x_0}\dd m
=c\int_0^1p_{x_0}\dd m=cx_0,
\]
a contradiction. Thus \(R\), and hence \(r\) almost everywhere,
is constant.
\end{proof}

\begin{proposition}\label{thm:strict}
Suppose that the strict triple \eqref{eq:convex-combination} has additive generators \(f_i\) such that \(f_i\) and \(f_i^{-1}\) are locally absolutely continuous in the interiors of their domains, for \(i=1,2,3\). Then
\[
T_1=T_2=T_3.
\]
\end{proposition}
\begin{proof}
In the nontrivial alternative of Lemma~\ref{lem:normalization}, Lemmas~\ref{lem:compatibility} and~\ref{lem:ratio} give a constant \(c>0\) with
\[
v_2=cv_1\ \text{a.e.},\qquad f_2'=f_1'/c\ \text{a.e.}
\]
For \(0<x<t<1\),
\(f_2(t)-f_2(x)=(f_1(t)-f_1(x))/c\). Let \(t\uparrow1\). Then \(f_2=f_1/c\), so \(T_2=T_1\) by generator uniqueness and \(T_3=T_1\) by \eqref{eq:convex-combination}. The other alternative of Lemma~\ref{lem:normalization} already gives the conclusion.
\end{proof}

\section{Common levels and a first change of coordinate}\label{sec:common}
Let \(P(u,v)=uv\). Fix strict triangular norms \(U,V\), an increasing homeomorphism \(\psi:[0,1]\to[0,1]\), and positive weights \(\lambda+\mu=1\), such that
\begin{equation}\label{eq:common-outer}
\psi(uv)=\lambda\psi(U(u,v))+\mu\psi(V(u,v)).
\end{equation}
Every strict triple satisfying \eqref{eq:convex-combination} gives this identity: choose a multiplicative generator \(h_3\) of \(T_3\), put \(\psi=h_3^{-1}\), and set
\begin{equation}\label{eq:original-conjugacy}
\begin{aligned}
U(u,v)&=h_3\bigl(T_1(\psi(u),\psi(v))\bigr),\\
V(u,v)&=h_3\bigl(T_2(\psi(u),\psi(v))\bigr).
\end{aligned}
\end{equation}
Neither the generators nor the outer coordinate \(\psi\) are assumed to be absolutely continuous.

Choose a multiplicative generator \(H\) of \(U\). For \(0<z<1\), define
\begin{equation}\label{eq:level-products}
n_z(v)=H^{-1}\!\left(\frac{H(z)}{H(v)}\right),\qquad
p_z(v)=v n_z(v),\qquad I_z=p_z([z,1]).
\end{equation}
For \(v\in[z,1]\),
\[
U(n_z(v),v)=z,\qquad
n_z(n_z(v))=v,\qquad n_z(z)=1,\quad n_z(1)=z.
\]
Thus \(n_z\) is a decreasing involution. Put
\begin{equation}\label{eq:common-level-set}
\begin{aligned}
S&=\{0,1\}\cup
 \{z\in(0,1):p_z(v)=z\text{ for all }v\in[z,1]\},\\
D&=(0,1)\setminus S.
\end{aligned}
\end{equation}
For \(0<z<1\), strictness gives
\[
z\in S\quad\Longleftrightarrow\quad
\{(u,v):U(u,v)=z\}=\{(u,v):uv=z\}.
\]
The map \((z,v)\mapsto p_z(v)\) is continuous on \(0<z\leq v\leq1\).
If \(z_n\in S\) and \(z_n\to z\in(0,1)\), then, for each fixed \(v>z\),
\(p_z(v)=\lim_n p_{z_n}(v)=\lim_nz_n=z\).
Also \(p_z(z)=z\), so \(z\in S\). Hence \(S\) is closed. No assumption on \(m(S)\) is made.

\begin{lemma}\label{lem:common-level-functions}
For every \(z\in S\) and \(u,v\in[0,1]\),
\begin{equation}\label{eq:common-thresholds}
uv\leq z\quad\Longleftrightarrow\quad U(u,v)\leq z
\quad\Longleftrightarrow\quad V(u,v)\leq z.
\end{equation}
If \(\xi\) is a finite nonnegative atomless Borel measure concentrated on \(S\), and \(b(t)=\xi([0,t])\), then
\begin{equation}\label{eq:common-cdf-values}
b(uv)=b(U(u,v))=b(V(u,v)).
\end{equation}
Consequently, adding \(b\) to \(\psi\), or subtracting it whenever the resulting function is nondecreasing, preserves the unnormalized identity \eqref{eq:common-outer}.
\end{lemma}
\begin{proof}
Fix \(z\in S\cap(0,1)\). If \(v\leq z\), all three outputs are at most \(v\).
If \(v>z\), put \(u_0=z/v\). Then
\[
U(u_0,v)=z,\qquad
\psi(z)=\lambda\psi(z)+\mu\psi(V(u_0,v)),
\]
so \(V(u_0,v)=z\). Strict increase of the sections gives
\[
\{u:uv\leq z\}
=\{u:U(u,v)\leq z\}
=\{u:V(u,v)\leq z\}=[0,u_0].
\]
For \(z=0\), each output is zero exactly when \(uv=0\); for \(z=1\), all outputs are at most~\(1\). This proves \eqref{eq:common-thresholds}.

Set \(a_0=uv,a_1=U(u,v),a_2=V(u,v)\). For \(j=0,1,2\), atomlessness and \(\xi([0,1]\setminus S)=0\) give
\[
b(a_j)=\xi([0,1])-\int_S\mathbf1_{\{a_j\leq z\}}\,\xi(\mathrm dz).
\]
The three integrands agree by \eqref{eq:common-thresholds}; hence \eqref{eq:common-cdf-values}. Since \(\lambda+\mu=1\),
\[
(\psi\pm b)(uv)
=\lambda(\psi\pm b)(U(u,v))+\mu(\psi\pm b)(V(u,v)).
\]
\end{proof}

\begin{lemma}\label{lem:involution-null}
Let \(0<c<1\), let \(n:[c,1]\to[c,1]\) be a decreasing involution, and put \(p(v)=v n(v)\). If a Borel set \(E\subseteq\mathbb R\) satisfies
\[
m\bigl(p^{-1}(E)\bigr)=0,
\]
then
\[
m\bigl(E\cap p([c,1])\bigr)=0.
\]
No absolute-continuity assumption is imposed on \(n\) or \(p\).
\end{lemma}
\begin{proof}
Since \(n^{-1}=n\), \(n\) is a homeomorphism. Put \(B=p^{-1}(E)\).
Then \(p\circ n=p\) and \(n(B)=B\).
Given \(\varepsilon>0\), choose a relatively open \(O\supset B\) with \(m(O)<\varepsilon\), and set
\[
W=O\cap n^{-1}(O),\qquad
B\subset W=n(W),\qquad m(W)<\varepsilon.
\]
Write \(W=\bigsqcup_k I_k\) as a countable union of interval components. Then
\(W=\bigsqcup_k n(I_k)\), and for \(v,w\in I_k\),
\[
|p(v)-p(w)|\leq |v-w|+|n(v)-n(w)|,
\qquad \operatorname{diam}p(I_k)\leq |I_k|+|n(I_k)|.
\]
Thus, with \(m^*\) denoting outer Lebesgue measure,
\[
m^*(p(B))\leq\sum_k\operatorname{diam}p(I_k)
\leq\sum_k\bigl(|I_k|+|n(I_k)|\bigr)
=2m(W)<2\varepsilon.
\]
Letting \(\varepsilon\downarrow0\) and using \(p(B)=E\cap p([c,1])\) proves the claim.
\end{proof}

\begin{lemma}\label{lem:localized-domination}
Let \(\rho=\mathrm d\psi\) in \eqref{eq:common-outer}. Then
\begin{equation}\label{eq:localized-domination}
m|_D\ll\rho.
\end{equation}
Equivalently, \(\rho(E)=0\) implies \(m(E\cap D)=0\) for every Borel set \(E\).
\end{lemma}
\begin{proof}
Fix a Borel set \(E\subset(0,1)\) with \(\rho(E)=0\). For \(0<v<1\), set
\(U_v(u)=U(u,v)\) and \(B_v=\{u\in(0,1):vu\in E\}\).
By \eqref{eq:common-outer} and \eqref{eq:stieltjes-image},
\[
0\leq\lambda\rho(U_v(B_v))
\leq\rho(vB_v)=\rho(E\cap(0,v))=0.
\]
Since \(vU_v^{-1}(z)=p_z(v)\),
\[
\begin{aligned}
0
&=\int_0^1\rho\{z\in(0,v):p_z(v)\in E\}\dd v\\
&=\int_{(0,1)}\left(\int_z^1\mathbf1_E(p_z(v))\dd v\right)\rho(\mathrm dz).
\end{aligned}
\]
Tonelli applies because \((z,v)\mapsto p_z(v)\) is continuous on \(0<z<v<1\). Hence
\[
G_E=\left\{z\in(0,1):\int_z^1\mathbf1_E(p_z(v))\dd v=0\right\}
\]
has full \(\rho\)-measure. The measure \(\rho\) has full support, so \(G_E\) is dense.
Adding the two parameter endpoints changes no Lebesgue measure. Lemma~\ref{lem:involution-null} therefore gives
\begin{equation}\label{eq:localized-null-interval}
m(E\cap I_z)=0\qquad(z\in G_E).
\end{equation}

Fix \(a\in D\). Choose \(v_0\in(a,1)\) with \(p_a(v_0)\ne a\), and put
\(\varepsilon=|p_a(v_0)-a|/2\). Continuity gives
\(0<\delta<\min\{a,v_0-a,1-a\}\) such that
\[
|z-a|<\delta\quad\Longrightarrow\quad
|p_z(v_0)-p_a(v_0)|<\varepsilon.
\]
If \(p_a(v_0)>a\), choose \(z\in G_E\cap(a-\delta,a)\); then
\(z<a<p_z(v_0)\). If \(p_a(v_0)<a\), choose
\(z\in G_E\cap(a,a+\delta)\); then \(p_z(v_0)<a<z\).
As \(p_z(1)=z\) and \(I_z\) is an interval, both cases give
\(a\in\operatorname{int}I_z\). Thus a countable subfamily satisfies
\[
D\subseteq\bigcup_{j\geq1}\operatorname{int}I_{z_j},
\qquad z_j\in G_E.
\]
By \eqref{eq:localized-null-interval},
\(m(E\cap D)\leq\sum_jm(E\cap I_{z_j})=0\).
\end{proof}

\begin{proposition}\label{prop:first-coordinate}
Define
\begin{equation}\label{eq:first-coordinate}
\ell(t)=m(S\cap[0,t]),\qquad
\psi_+(t)=\frac{\psi(t)+\ell(t)}{1+m(S)}.
\end{equation}
Then \(\psi_+\) is an increasing homeomorphism of \([0,1]\),
\begin{equation}\label{eq:first-coordinate-combination}
\psi_+(uv)=\lambda\psi_+(U(u,v))+\mu\psi_+(V(u,v)),
\end{equation}
and \(m\ll\mathrm d\psi_+\).
\end{proposition}
\begin{proof}
The function \(\ell\) is nondecreasing and Lipschitz, so \(\psi_+\) is a continuous strictly increasing map with endpoints \(0,1\). Lemma~\ref{lem:common-level-functions}, applied to \(m|_S\), gives
\(\ell(uv)=\ell(U(u,v))=\ell(V(u,v))\), proving \eqref{eq:first-coordinate-combination}. Moreover,
\[
\rho_+=\mathrm d\psi_+=\frac{\rho+\mathbf1_Sm}{1+m(S)}.
\]
For a Borel set \(E\), Lemma~\ref{lem:localized-domination} yields
\[
\rho_+(E)=0
\ \Longrightarrow\ \rho(E)=m(E\cap S)=0
\ \Longrightarrow\ m(E\cap D)=0
\ \Longrightarrow\ m(E)=0.
\]
Thus \(m\ll\rho_+\).
\end{proof}

\section{Absolute continuity of the relative generators}\label{sec:forced}
We work directly with \eqref{eq:common-outer}. Put \(U_1=U\), \(U_2=V\), \(U_3=P\), and \(\rho=\mathrm d\psi\). Choose multiplicative generators \(H_i\) of \(U_i\), \(i=1,2\). For \(v\in(0,1)\), define the section measures
\begin{equation}\label{eq:kappa-definition}
\kappa_{i,v}=\mathrm d_u[\psi(U_i(u,v))]\qquad(i=1,2,3).
\end{equation}
For every Borel set \(B\subseteq[0,1]\),
\begin{equation}\label{eq:kappa-images}
\kappa_{i,v}(B)=\rho(U_i(B,v)),\qquad
\kappa_{3,v}(B)=\rho(vB),\qquad
\kappa_{3,v}=\lambda\kappa_{1,v}+\mu\kappa_{2,v}.
\end{equation}
Since \(U_i(\cdot,v)\to\operatorname{id}\) uniformly as \(v\uparrow1\), every compact \(K\subset(0,1)\) satisfies
\begin{equation}\label{eq:compact-upper}
\limsup_{v\uparrow1}\kappa_{i,v}(K)\leq\rho(K).
\end{equation}
Indeed, for each \(\varepsilon>0\), the image \(U_i(K,v)\) eventually lies in the closed \(\varepsilon\)-neighbourhood \(K_\varepsilon\) of \(K\), and \(\rho(K_\varepsilon)\downarrow\rho(K)\).

We use the standard continuity of dilations in \(L^1(0,\infty)\):
\begin{equation}\label{eq:dilation-L1}
\|a g(a\,\cdot)-b g(b\,\cdot)\|_1\longrightarrow0
\qquad(a\to b>0).
\end{equation}
Under the isometry \(G(s)=e^s g(e^s)\), the norm on the left equals
\(\|G(\cdot+\log a)-G(\cdot+\log b)\|_{L^1(\mathbb R)}\), so this is translation continuity \cite[Proposition~8.5]{Folland1999}. Functions on \((0,1)\) are extended by zero when this fact is used.

The following two consequences of Tonelli's theorem provide the required section estimates.

\begin{lemma}\label{lem:dilation}
Let \(\eta\) be a finite positive Borel measure on \((0,1)\).
\begin{enumerate}
\item If \(N\subset(0,1)\) is a Lebesgue-null Borel set, then \(\eta(cN)=0\) for almost every \(c\in(0,1)\).
\item If \(\eta\) has full support, \(k\in L^1(0,1)\) is nonnegative and nonzero, and \(B\subset(0,1)\) is compact with \(m(B)>0\), then
\begin{equation}\label{eq:positive-average}
\int_0^1\eta(cB)k(c)\dd c>0.
\end{equation}
\end{enumerate}
\end{lemma}
\begin{proof}
For~(i), fix \(z\in(0,1)\). The set
\(\{c\in(0,1):z/c\in N\}=\{z/t:t\in N\cap(z,1)\}\)
is null, since \(t\mapsto z/t\) is Lipschitz on \([z,1]\). Tonelli's theorem gives
\[
\int_0^1\eta(cN)\dd c
=\int_{(0,1)}\left(\int_0^1\mathbf1_N(z/c)\dd c\right)\eta(\mathrm dz)=0.
\]

For~(ii), extend \(k\) by zero and set
\[
G(z)=\int_0^\infty k(c)\mathbf1_B(z/c)\dd c
=\int_0^\infty z k(zt)\mathbf1_B(1/t)\dd t.
\]
Equation~\eqref{eq:dilation-L1} makes \(G\) continuous for \(z>0\). Tonelli's theorem also gives
\[
\int_0^1G(z)\dd z=m(B)\int_0^1c k(c)\dd c>0,
\qquad
\int_0^1\eta(cB)k(c)\dd c=\int_{(0,1)}G(z)\,\eta(\mathrm dz).
\]
Thus \(G\) is bounded below by a positive constant on some nonempty open interval. Full support of \(\eta\) makes the last integral positive.
\end{proof}

\begin{proposition}\label{prop:relative-regularity-direct}
Suppose that \eqref{eq:common-outer} holds and
\(m\ll\rho=\mathrm d\psi\). If \(H_1,H_2\) are multiplicative
generators of \(U,V\), respectively, then
\[
H_i,H_i^{-1}\in\AC([0,1])\qquad(i=1,2).
\]
\end{proposition}
\begin{proof}
Fix \(i\in\{1,2\}\), put \(j=3-i\), \(H=H_i\), and
\((w_1,w_2)=(\lambda,\mu)\).
Write \(\rho=r\,m+\nu\) and choose a Borel null set \(N\)
carrying \(\nu\). Since
\(\rho(\{r=0\}\setminus N)=0\) and \(m\ll\rho\),
we have \(r>0\) almost everywhere. Set
\(\eta=\mathrm d(\psi\circ H^{-1})\); this is a finite
atomless measure with full support.

\emph{Absolute continuity of \(H^{-1}\).}
Suppose otherwise. Choose a null Borel set \(A\subset(0,1)\)
with \(m(H^{-1}(A))>0\). Inner regularity gives a compact set
\[
K\subset H^{-1}(A)\setminus N,\qquad
0<\rho(K)=\int_Kr\,\mathrm dm,
\qquad m(H(K))=0.
\]
Lemma~\ref{lem:dilation}(i) supplies \(c_n\uparrow1\) with
\(\eta(c_nH(K))=0\). Put \(v_n=H^{-1}(c_n)\uparrow1\).
Then
\[
\kappa_{i,v_n}(K)=\eta(c_nH(K))=0,
\qquad \kappa_{3,v_n}(K)=w_j\kappa_{j,v_n}(K).
\]
On the other hand,
\[
\kappa_{3,v}(K)=\rho(vK)\geq\int_Kv r(vu)\,\mathrm du,
\qquad
\left|\int_Kv r(vu)\,\mathrm du-\rho(K)\right|
\leq\|v r(v\,\cdot)-r\|_1\longrightarrow0.
\]
Together with \eqref{eq:compact-upper}, this yields
\[
\rho(K)\leq\liminf_n\kappa_{3,v_n}(K)
\leq w_j\limsup_n\kappa_{j,v_n}(K)
\leq w_j\rho(K)<\rho(K),
\]
a contradiction. Hence \(H^{-1}\in\AC([0,1])\).

\emph{Absolute continuity of \(H\).}
If \(H\notin\AC\), inner regularity gives a compact
\(K\subset(0,1)\) with \(m(K)=0<m(H(K))\).
Lemma~\ref{lem:dilation}(i) and
\(w_i\kappa_{i,v}\leq\kappa_{3,v}\) give
\[
\rho(vK)=0\quad\text{for almost every }v,
\qquad \int_0^1\kappa_{i,v}(K)\,\mathrm dv=0.
\]
Put \(B=H(K)\) and \(k=(H^{-1})'\). Then
\(k\geq0\), \(\int_0^1k=1\), and the pushforward of
\(m\) under \(H\) is \(k(c)\,\mathrm dc\), because
\[
m\{v:H(v)\in(a,b]\}
=H^{-1}(b)-H^{-1}(a)=\int_a^bk(c)\,\mathrm dc.
\]
Consequently
\[
0=\int_0^1\eta(H(v)B)\,\mathrm dv
=\int_0^1\eta(cB)k(c)\,\mathrm dc>0,
\]
where the strict inequality is Lemma~\ref{lem:dilation}(ii).
Thus \(H\in\AC([0,1])\). The argument applies to both indices.
\end{proof}

\section{Combination of strict triangular norms}\label{sec:regularization}
Return to \eqref{eq:common-outer} and the closed set \(S\) from Section~\ref{sec:common}. The next lemma shows that absolute continuity of a relative generator and its inverse concentrates the singular part of \(\mathrm d\psi\) on \(S\). We then apply it to \(\psi_+\).

\begin{lemma}\label{lem:singular-support}
Suppose \eqref{eq:common-outer} holds and a generator \(H\) of \(U\) and its inverse are absolutely continuous. Write
\[
\psi=q+s,\qquad q(t)=\int_0^t\psi'(r)\dd r,
\]
where \(s\) is continuous, nondecreasing and singular, and put \(\nu=\mathrm ds\). Then \(\nu((0,1)\setminus S)=0\). In particular,
\begin{equation}\label{eq:q-combination}
q(uv)=\lambda q(U(u,v))+\mu q(V(u,v)).
\end{equation}
\end{lemma}
\begin{proof}
Choose a Borel set \(N\subset(0,1)\) with
\(m(N)=0\) and \(\nu((0,1)\setminus N)=0\).
For \(0<v<1\), put \(U_v(u)=U(u,v)\).
The generator formula and Proposition~\ref{prop:ACfacts} give
\(U_v^{-1}\in\AC\). Hence \(\mathrm d(s\circ U_v)\)
is singular, being concentrated on the null set
\(U_v^{-1}(N\cap[0,v])\).
Since
\[
0\leq\lambda\,\mathrm d(s\circ U_v)
\leq\lambda\,\mathrm d(\psi\circ U_v)
\leq\mathrm d(\psi(v\,\cdot)),
\]
taking singular parts gives
\[
\lambda\,\mathrm d(s\circ U_v)
\leq\mathrm d(s(v\,\cdot)).
\]
The right-hand measure is concentrated on
\(E_v=\{u:vu\in N\}\). Using
\(vU_v^{-1}(z)=p_z(v)\), we obtain
\[
0=\mathrm d(s\circ U_v)(E_v^c)
=\nu\{z\in(0,v):p_z(v)\notin N\}.
\]
Tonelli gives
\[
\int_{(0,1)}m\{v\in(z,1):p_z(v)\notin N\}
\,\nu(\mathrm dz)=0.
\]
Thus, for \(\nu\)-almost every \(z\), the inverse image
\(p_z^{-1}([0,1]\setminus N)\subset[z,1]\) is null;
including the two endpoints does not change its measure.
Lemma~\ref{lem:involution-null} yields
\[
m(I_z\setminus N)=0,
\qquad m(I_z)=0.
\]
Since \(I_z\) is an interval containing \(z\), it equals
\(\{z\}\), so \(z\in S\). Hence
\(\nu((0,1)\setminus S)=0\).

Finally, \(s(0)=0\) and atomlessness give
\(s(t)=\nu([0,t])\). Lemma~\ref{lem:common-level-functions}
yields \(s(uv)=s(U(u,v))=s(V(u,v))\).
Subtracting these equal values from \eqref{eq:common-outer}
proves \eqref{eq:q-combination}.
\end{proof}

\begin{theorem}\label{thm:regularization}
Assume \eqref{eq:common-outer}. There is an increasing homeomorphism \(Q:[0,1]\to[0,1]\), with \(Q,Q^{-1}\in\AC([0,1])\), such that
\begin{equation}\label{eq:final-combination}
Q(uv)=\lambda Q(U(u,v))+\mu Q(V(u,v)).
\end{equation}
The multiplicative generators \(H_1,H_2\) of \(U,V\) and their inverses are also absolutely continuous.
\end{theorem}
\begin{proof}
Proposition~\ref{prop:first-coordinate} gives \(\psi_+\) satisfying \eqref{eq:first-coordinate-combination} and \(m\ll\rho_+=\mathrm d\psi_+\). Proposition~\ref{prop:relative-regularity-direct} therefore gives \(H_i,H_i^{-1}\in\AC\), \(i=1,2\).

Write \(\rho_+=r\,m+\nu\), put \(q_+(t)=\int_0^tr(u)\dd u\), and choose a Borel null set \(N\) carrying \(\nu\). Since
\(\rho_+(\{r=0\}\setminus N)=0\) and \(m\ll\rho_+\), we have \(r>0\) almost everywhere. Hence
\[
c=q_+(1)>0,\qquad Q(t)=q_+(t)/c
\]
defines an absolutely continuous increasing homeomorphism with \(Q'>0\) almost everywhere. Proposition~\ref{prop:ACfacts}(ii) gives \(Q^{-1}\in\AC\). Lemma~\ref{lem:singular-support}, applied to \(\psi_+\), shows that \(q_+\) satisfies the convex identity. Dividing by \(c\) proves \eqref{eq:final-combination}.
\end{proof}

The construction gives the explicit formula
\begin{equation}\label{eq:explicit-Q}
Q(t)=\frac{\displaystyle\int_0^t\psi'(r)\dd r+m(S\cap[0,t])}
{\displaystyle\int_0^1\psi'(r)\dd r+m(S)},
\end{equation}
whose denominator is positive by the proof.

\begin{proof}[Proof of Theorem~\ref{thm:full-strict}]
Only the forward implication needs proof. Suppose that
\(T_3=\lambda T_1+\mu T_2\) is associative, where \(\mu=1-\lambda\). It is a strict triangular norm. Choose its multiplicative generator \(h_3\), and define \(\psi,U,V\) by \eqref{eq:original-conjugacy}.

Take \(Q,H_1,H_2\) from Theorem~\ref{thm:regularization}, and set \(U_1=U\), \(U_2=V\), \(U_3=P\), \(H_3=\operatorname{id}\). The strict triangular norms
\[
\widehat T_i(x,y)=Q\bigl(U_i(Q^{-1}(x),Q^{-1}(y))\bigr)
\qquad(i=1,2,3)
\]
satisfy \(\widehat T_3=\lambda\widehat T_1+\mu\widehat T_2\) by \eqref{eq:final-combination}. Their generators and inverse generators are
\[
g_i=H_i\circ Q^{-1},\qquad g_i^{-1}=Q\circ H_i^{-1},
\]
which are absolutely continuous by Proposition~\ref{prop:ACfacts}(iii). Thus \(f_i=-\log g_i\) and \(f_i^{-1}(s)=g_i^{-1}(e^{-s})\) satisfy \eqref{eq:genregular}. Proposition~\ref{thm:strict} gives \(\widehat T_1=\widehat T_2=\widehat T_3\). Since \(Q\) is injective, \(U=V=P\); conjugating back by \(h_3\) gives \(T_1=T_2=T_3\).
\end{proof}

\section{Consequences}\label{sec:scope}
The ordinal-sum reduction of Ouyang and Fang \cite{OF2008} and the nilpotent result of Petr\'ik and Sarkoci \cite{PS2009} extend Theorem~\ref{thm:full-strict} to all continuous triangular norms.

\begin{corollary}\label{cor:continuous}
Let \(T_1,T_2\) be continuous triangular norms and let \(0<\lambda<1\). Then
\[
T_\lambda=\lambda T_1+(1-\lambda)T_2
\]
is associative if and only if \(T_1=T_2\).
\end{corollary}
\begin{proof}
The converse is immediate. Suppose that \(T_3=T_\lambda\) is associative and put \(\mu=1-\lambda\). Since the other triangular-norm axioms and continuity are preserved by convex combinations, \(T_3\) is a continuous triangular norm. Assume that \(T_1\ne T_2\).

For a continuous triangular norm \(A\), write
\[
E(A)=\{t\in[0,1]:A(t,t)=t\}.
\]
The ordinal-sum reduction of \cite{OF2008}, as stated in \cite[p.~2 of the expanded version]{NPS2009}, gives
\(E(T_1)=E(T_2)=E\). Since \(T_i(t,t)\leq t\) and \(\lambda,\mu>0\),
\[
\begin{aligned}
t\in E(T_3)
&\Longleftrightarrow \lambda T_1(t,t)+\mu T_2(t,t)=t\\
&\Longleftrightarrow T_1(t,t)=T_2(t,t)=t.
\end{aligned}
\]
Hence \(E(T_3)=E\).

By the ordinal-sum representation of continuous triangular norms \cite{KMP2000}, each \(T_i\) equals \(\min\) outside the squares \([a,b]^2\) associated with the components \((a,b)\) of \([0,1]\setminus E\). Since \(T_1\ne T_2\), there is such a component on whose square they differ. Put \(d=b-a\), and define
\[
S_i(u,v)=\frac{T_i(a+du,a+dv)-a}{d},
\qquad (u,v)\in[0,1]^2,\quad i=1,2,3.
\]
These are continuous Archimedean triangular norms, \(S_1\ne S_2\), and
\[
S_3=\lambda S_1+\mu S_2.
\]
If both \(S_1,S_2\) are nilpotent, this contradicts \cite{PS2009}. If exactly one is nilpotent, it contradicts the mixed-case obstruction of \cite{OF2008}; see also \cite[Theorems~1 and~2]{NPS2009}. The classification of continuous Archimedean triangular norms \cite{KMP2000} therefore implies that both \(S_1,S_2\) are strict. Theorem~\ref{thm:full-strict} gives \(S_1=S_2\), a contradiction.
\end{proof}

Corollary~\ref{cor:continuous} answers \cite[Problem~1]{Navara2015} and the continuous version of \cite[Problem~5]{AFS2003}. Its proof uses a single component on which the operations differ and therefore applies to infinite ordinal sums as well. Common nilpotent components are allowed, but the two operations must agree on each such component.

A \emph{bivariate copula} is a map \(C:[0,1]^2\to[0,1]\) satisfying
\[
C(u,0)=C(0,u)=0,\qquad C(u,1)=C(1,u)=u
\]
for \(u\in[0,1]\), and
\[
C(u_2,v_2)-C(u_1,v_2)-C(u_2,v_1)+C(u_1,v_1)\geq0
\]
whenever \(0\leq u_1\leq u_2\leq1\) and \(0\leq v_1\leq v_2\leq1\); see \cite[Definition~5]{KMPGenerated2000}. These conditions are preserved by convex combinations.

\begin{corollary}\label{cor:copulas}
Let \(C_1,C_2\) be associative bivariate copulas and let \(0<\lambda<1\). Then
\[
C_\lambda=\lambda C_1+(1-\lambda)C_2
\]
is associative if and only if \(C_1=C_2\).
\end{corollary}
\begin{proof}
Every copula is continuous, and every associative copula is commutative \cite[pp.~370--371, Proposition~3]{KMPGenerated2000}. Thus \(C_1,C_2\) are continuous triangular norms, and Corollary~\ref{cor:continuous} applies.
\end{proof}

Corollary~\ref{cor:copulas} gives a negative answer to \cite[Open Question~1, p.~237]{Seliga2026}. In particular, for the product copula \(\TP(x,y)=xy\),
\[
\TP=\lambda C_1+(1-\lambda)C_2
\quad\Longrightarrow\quad C_1=C_2=\TP
\]
under the hypotheses of the corollary. This answers the two-summand product-decomposition question in \cite[Problem~6]{AFS2003}.

A \emph{triangular conorm} is a commutative, associative and nondecreasing operation on \([0,1]\) with neutral element \(0\). Its dual \(T(x,y)=1-S(1-x,1-y)\) is a triangular norm \cite{KMP2000}.

\begin{corollary}\label{cor:conorm}
Let \(S_1,S_2\) be continuous triangular conorms and let \(0<\lambda<1\). Then \(\lambda S_1+(1-\lambda)S_2\) is associative if and only if \(S_1=S_2\).
\end{corollary}
\begin{proof}
The duality \(\mathcal D F(x,y)=1-F(1-x,1-y)\) preserves associativity and continuity and commutes with convex combinations \cite[Remark~2]{Navara2015}. Apply Corollary~\ref{cor:continuous} to \(\mathcal D S_1,\mathcal D S_2\).
\end{proof}

\begin{corollary}\label{cor:petrik-full}
Petr\'ik's Conjecture~1 holds. Its assumptions on boundary derivatives and on the absolute continuity of the canonical generators are unnecessary within the class of strict triangular norms.
\end{corollary}
\begin{proof}
Conjecture~1 of \cite[p.~1055]{Petrik2010} is a restricted instance of Theorem~\ref{thm:full-strict}.
\end{proof}

A \emph{uninorm} is a commutative, associative and nondecreasing operation on \([0,1]\) with a neutral element \(e\in[0,1]\). For \(0<e<1\), it is \emph{representable} if there is a continuous strictly increasing bijection \(g:(0,1)\to\mathbb R\), with \(g(e)=0\), such that
\[
U(x,y)=g^{-1}(g(x)+g(y))\qquad(0<x,y<1),
\]
and
\[
U(x,0)=0\quad(0\leq x<1),\qquad
U(x,1)=1\quad(0<x\leq1),\qquad
U(0,1)=U(1,0)\in\{0,1\};
\]
see \cite[Section~2]{MZ2021}. The following consequence uses the reduction in \cite[Proposition~33]{MZ2021}.

\begin{corollary}\label{cor:representable-uninorms}
Let \(U_1,U_2\) be representable uninorms with the same neutral element \(e\in(0,1)\), and let \(0<\lambda<1\). Then \(\lambda U_1+(1-\lambda)U_2\) is a uninorm if and only if \(U_1=U_2\).
\end{corollary}
\begin{proof}
Proposition~33 of \cite{MZ2021} proves this assertion under the hypothesis that no convex combination with positive coefficients of two distinct strict triangular norms is a triangular norm. That hypothesis is precisely Theorem~\ref{thm:full-strict}. The converse is immediate.
\end{proof}

The common-neutral-element hypothesis is the one used in \cite[Proposition~33]{MZ2021}; the questions with different neutral elements in \cite[Section~7]{MZ2021} remain separate.

The following standard construction illustrates why the proof regularizes an auxiliary triple; compare \cite[Exercise~41(a), p.~109]{Folland1999}.

\begin{example}\label{ex:degenerate-generator}
Let \(K\subset(1/3,2/3)\) be compact and nowhere dense, with \(m(K)>0\), and put
\[
h(x)=\frac{1}{1-m(K)}\int_0^x\mathbf1_{[0,1]\setminus K}(t)\dd t.
\]
Every nonempty interval meets the open complement of \(K\), so \(h\) is an absolutely continuous increasing homeomorphism of \([0,1]\). Since \(h'=0\) almost everywhere on \(K\), Proposition~\ref{prop:ACfacts}(ii) gives \(h^{-1}\notin\AC\). The strict norm generated by \(h\) is nevertheless covered by Theorem~\ref{thm:full-strict}.
\end{example}

Continuity of both summands in Corollary~\ref{cor:continuous} is essential.

\begin{example}\label{ex:discontinuous}
Let \(P(x,y)=xy\) and let \(D\) be the drastic product. The classical family \cite[Example~1]{Navara2015}
\[
F_\lambda(x,y)=\lambda P(x,y)+(1-\lambda)D(x,y)
=\begin{cases}
\lambda xy,&x,y<1,\\
\min\{x,y\},&\max\{x,y\}=1,
\end{cases}
\]
consists of triangular norms for \(0<\lambda<1\). Each is discontinuous, since \(F_\lambda(x,y)\to\lambda x\ne x=F_\lambda(x,1)\) as \(y\uparrow1\), for \(0<x<1\).
\end{example}

\setlength{\bibsep}{4pt}


\begin{thebibliography}{99}

\bibitem{KMP2000}
E.P. Klement, R. Mesiar, E. Pap,
\emph{Triangular Norms}, Trends in Logic, vol.~8,
Kluwer Academic Publishers, Dordrecht, 2000.
\href{https://doi.org/10.1007/978-94-015-9540-7}{doi:10.1007/978-94-015-9540-7}.

\bibitem{KMPGenerated2000}
E.P. Klement, R. Mesiar, E. Pap,
Generated triangular norms,
Kybernetika 36 (2000) 363--377.
\href{https://www.kybernetika.cz/content/2000/3/363/paper.pdf}{Full text}.

\bibitem{AFS2003}
C. Alsina, M.J. Frank, B. Schweizer,
Problems on associative functions,
Aequationes Math. 66 (2003) 128--140.
\href{https://doi.org/10.1007/s00010-003-2673-y}{doi:10.1007/s00010-003-2673-y}.

\bibitem{Jenei2006}
S. Jenei,
On the convex combination of left-continuous \(t\)-norms,
Aequationes Math. 72 (2006) 47--59.
\href{https://doi.org/10.1007/s00010-006-2840-z}{doi:10.1007/s00010-006-2840-z}.

\bibitem{DS2008}
F. Durante, P. Sarkoci,
A note on the convex combinations of triangular norms,
Fuzzy Sets Syst. 159 (2008) 77--80.
\href{https://doi.org/10.1016/j.fss.2007.07.005}{doi:10.1016/j.fss.2007.07.005}.

\bibitem{MM2008}
R. Mesiar, A. Mesiarov\'a-Zem\'ankov\'a,
Convex combinations of continuous \(t\)-norms with the same diagonal function,
Nonlinear Anal. 69 (2008) 2851--2856.
\href{https://doi.org/10.1016/j.na.2007.08.076}{doi:10.1016/j.na.2007.08.076}.

\bibitem{OF2008}
Y. Ouyang, J. Fang,
Some observations about the convex combination of continuous triangular norms,
Nonlinear Anal. 68 (2008) 3382--3387.
\href{https://doi.org/10.1016/j.na.2007.03.027}{doi:10.1016/j.na.2007.03.027}.

\bibitem{PS2009}
M. Petr\'ik, P. Sarkoci,
Convex combinations of nilpotent triangular norms,
J. Math. Anal. Appl. 350 (2009) 271--275.
\href{https://doi.org/10.1016/j.jmaa.2008.09.060}{doi:10.1016/j.jmaa.2008.09.060}.

\bibitem{Petrik2010}
M. Petr\'ik,
Convex combinations of strict \(t\)-norms,
Soft Comput. 14 (2010) 1053--1057.
\href{https://doi.org/10.1007/s00500-009-0484-3}{doi:10.1007/s00500-009-0484-3}.

\bibitem{NPS2009}
M. Navara, M. Petr\'ik, P. Sarkoci,
Convex combinations of triangular norms,
in: U. Bodenhofer, B. De Baets, E.P. Klement, S. Saminger-Platz (Eds.),
\emph{30th Linz Seminar on Fuzzy Set Theory: The Legacy of 30 Seminars---Where Do We Stand and Where Do We Go?},
Johannes Kepler University, Linz, 2009, pp.~85--87.
Expanded author version: \href{https://www.flll.jku.at/div/research/linz2009/wwthmga_09_stcartsba/NavaraPetrikSarkoci_linz09.pdf}{four-page manuscript}.

\bibitem{Navara2015}
M. Navara,
Convex combinations of fuzzy logical operations,
Fuzzy Sets Syst. 264 (2015) 51--63.
\href{https://doi.org/10.1016/j.fss.2014.10.013}{doi:10.1016/j.fss.2014.10.013}.

\bibitem{Navara2019}
M. Navara,
Formulas for generators of \(R\)-implications,
Fuzzy Sets Syst. 359 (2019) 80--89.
\href{https://doi.org/10.1016/j.fss.2018.09.011}{doi:10.1016/j.fss.2018.09.011}.

\bibitem{NP2020}
M. Navara, M. Petr\'ik,
Generators of fuzzy logical operations,
in: H.T. Nguyen, V. Kreinovich (Eds.),
\emph{Algebraic Techniques and Their Use in Describing and Processing Uncertainty},
Studies in Computational Intelligence, vol.~878,
Springer, Cham, 2020, pp.~89--112.
\href{https://doi.org/10.1007/978-3-030-38565-1_8}{doi:10.1007/978-3-030-38565-1\_8}.

\bibitem{MZ2021}
A. Mesiarov\'a-Zem\'ankov\'a,
Convex combinations of uninorms and triangular subnorms,
Fuzzy Sets Syst. 423 (2021) 55--73.
\href{https://doi.org/10.1016/j.fss.2020.10.011}{doi:10.1016/j.fss.2020.10.011}.

\bibitem{Seliga2026}
A. \v{S}eliga,
A Schur-concave shuffle of \(M\) outside the closed convex hull of associative copulas,
Iranian J. Fuzzy Syst. 23 (6) (2026) 223--238.
\href{https://doi.org/10.22111/ijfs.2026.55340.9816}{doi:10.22111/ijfs.2026.55340.9816}.

\bibitem{Folland1999}
G.B. Folland,
\emph{Real Analysis: Modern Techniques and Their Applications},
second ed., John Wiley \& Sons, New York, 1999.

\bibitem{Leoni2009}
G. Leoni,
\emph{A First Course in Sobolev Spaces},
Graduate Studies in Mathematics, vol.~105,
American Mathematical Society, Providence, RI, 2009.
\href{https://doi.org/10.1090/gsm/105}{doi:10.1090/gsm/105}.

\end{thebibliography}
\end{document}